\documentclass[a4paper, 12pt] {article}
\usepackage[cp1251]{inputenc}
\usepackage[russian, english]{babel}
\usepackage{amsfonts}
\usepackage{mathtext}
\usepackage{cite}

\usepackage{graphicx}

\usepackage{verbatim}
\usepackage{amsmath}
\usepackage{amsthm}
\usepackage{amssymb}
\usepackage{delarray}
\usepackage{mathrsfs}
\usepackage{xcolor}

\begin{document}

\thispagestyle{empty}

\begin{center}
{\Large\bf Functional-differential operators on geometrical graphs with global delay and inverse spectral problems }

\end{center}

\begin{center}
{\large\bf Sergey Buterin\footnote{Department of Mathematics, Saratov State University, Russia {\it email: buterinsa@sgu.ru}} }
\end{center}

{\bf Abstract.} We suggest a new concept of functional-differential operators with constant delay on geometrical graphs that involves {\it
global} delay parameter. Differential operators on graphs model various processes in many areas of science and technology.
Although a vast majority of studies in this direction concern purely differential operators on graphs (often referred to as quantum graphs),
recently there also appeared some considerations of nonlocal operators on star-type graphs. In particular, there
belong functional-differential operators with constant delays but in a {\it locally} nonlocal version.
The latter means that each edge of the graph has its own delay parameter, which does not affect any other edge.
In this paper, we introduce {\it globally} nonlocal operators that are expected to be more natural for modelling
nonlocal processes on graphs. We also extend this idea to arbitrary trees, which opens a wide area of further research.
Another goal of the paper is to study inverse spectral problems for operators with global delay in one illustrative case by addressing
a wide range of questions including uniqueness, characterization of the spectral data as well as the uniform stability.

\smallskip
Key words: functional-differential equation, constant delay, globally nonlocal operator, metric graph, quantum graph,
inverse spectral problem, characterization, uniform stability

\smallskip
2010 Mathematics Subject Classification: 34A55 34K29 34B45\\
\\

{\large\bf 1. Introduction}
\\

Differential operators on geometrical graphs (spacial networks)
often appear in mathematics, mechanics, physics, organic chemistry, nanotechnology and other fields
of science and engineering (see \cite{Mont, Nic, vB, LangLeug, Kuch, BCFK, BerkKuch, Pok, Kuz-17}
and references therein). Metric graphs equipped with differential equations on
their edges are frequently referred to as {\it quantum graphs}, which indicates applications in quantum mechanics but
became common for such objects. In recent decades, theory of quantum graphs was significantly supplemented by
studies of inverse spectral problems (see, e.g., \cite{Ger, Bel, BrauWeik, Yur-05, Yurko-16, Bond18, Bond-20}).

Meanwhile, there were only few studies of {\it functional}-differential operators as well as other classes of nonlocal operators on graphs.
The definition
of such operators is complicated by some obvious difficulties caused by behavior of the nonlocality at internal vertices of the graph.
That is why the existing studies mostly address only {\it locally} nonlocal case when the corresponding nonlocal equation on each edge
can be considered independently on the other edges \cite{Nizh-12, Bon18-1, HuBondShYan19, Hu20, WangYang-21, Bon22, WangYang-22}.

The present paper is aimed to suggest another concept of functional-differential operators with delay on geometrical graphs
that involves global delay parameter. We believe that such settings might give ways to modelling natural nonlocal processes on various
branching structures. Although our idea can be extended in one or another way to many classes of nonlocal operators, we focus here
on Sturm--Liouville-type operators with constant delay of the form
\begin{equation}\label{1.1}
-y''(x)+q(x)y(x-a).
\end{equation}
Operators with delay as well as other classes of operators with deviating argument have been
actively studied starting from the middle of the
last century in connection with numerous applications \cite{Mysh, BellCook, Nor, Hal, Skub, Mur2}.
Recently, there appeared also many studies devoted to various aspects of inverse
problems for operator (1) on an interval \cite{Pik91, FrYur12, Yang, Ign18, BondYur18-1, ButYur19, VPV19, DV19, SatShieh19, WShM19, Dur20,
VPVC20, DB21, DB21-2, DB22, BMSh21}. Another our aim is to study inverse spectral problems for the introduced globally nonlocal
operators on one simple but illustrative graph, for which we intend to cover a wide range of specific questions.

\medskip
First, let us recall how the classical Sturm--Liouville operator (i.e. when $a=0)$ can be defined on a graph.
For this purpose, it is sufficient to consider star graph $G_m$ (see Fig.~1).
\begin{center}
\unitlength=1.0mm
\begin{picture}(80,55)
 \put(40,25){\line(-2,-3){15}}
 \put(40,25){\line(-5,2){25}}
 \put(40,25){\line(0,1){27}}
 \multiput(40,25)(2.5,1){11}{\circle*{0.7}}
 \put(40,25){\line(2,-3){15}}

 \put(40,25){\circle*{1}}
 \put(24.9,2.2){\circle*{1}}
 \put(14.7,35.1){\circle*{1}}
 \put(39.97,52){\circle*{1}}
 \put(55.1,2.2){\circle*{1}}

 \put (32.5,10){$e_1$}
 \put (23,26.5){$e_2$}
 \put (34.5,39.5){$e_3$}
 \put (50,13){$e_m$}

 \put (33,22.7){$v_0$}
 \put (20.5,-1.5){$v_1$}
 \put (9,35.5){$v_2$}
 \put (38.5,54){$v_3$}
 \put (56,-1.5){$v_m$}

  \put (70,0){\small Fig. 1. Graph $G_m$}
\end{picture}
\end{center}
Let each edge $e_j$ of $G_m$ be parameterized by the variable $x\in[0,1]$ so that $x=1$ always corresponds to the unique internal
vertex $v_0,$ while $x=0$ is associated with the corresponding boundary one $v_j.$ For $j=\overline{1,m},$ any function $f(x)$ on the edge
$e_j$ is actually defined on the segment $[0,1]$ and will be marked with the same index $j$ as, e.g., $f_j(x).$

Consider the boundary value problem on $G_m$ consisting of the Sturm--Liouville equations
\begin{equation}\label{1.1-0}
-y_j''(x)+q_j(x)y_j(x)=\lambda y_j(x), \quad 0<x<1, \quad j=\overline{1,m},
\end{equation}
on all its edges along with $m$ matching conditions at the internal vertex $v_0:$
\begin{equation}\label{1.3}
y_1(1)=y_2(1)=\ldots= y_m(1), \quad \sum_{j=1}^m y_j'(1)=0
\end{equation}
and $m$ (e.g., Dirichlet) boundary conditions at the boundary vertices:
\begin{equation}\label{1.4}
y_j(0)=0, \quad j=\overline{1,m},
\end{equation}
where $\lambda$ is the spectral parameter, while $q_j(x)$ are potentials in an appropriate class.

The specificity of considering differential equations on graphs is related to how their solutions match at internal vertices.
Conditions (\ref{1.3}), consisting of $m-1$ continuity conditions along with one Kirchhoff's condition
and often referred to as standard matching conditions, appear in many applications.
For example, they express the balance of tension in $m$ connected strings.

\medskip
Before discussing possible definitions of the operator (\ref{1.1}) with a nonzero value of $a$ on graphs, we point
out some specifics of the corresponding equation on an interval
\begin{equation}\label{1.1-1}
-y''(x)+q(x)y(x-a)=\lambda y(x), \quad 0<x<1,
\end{equation}
where $a\in(0,1)$ is fixed, while $q(x)$ is an integrable function.
Unlike the local case $a=0,$ equation (\ref{1.1-1}) is underdetermined since the argument of the unknown function $y$ may
fall out of the interval. There are two general ways how to overcome this issue, namely:

(i) To specify an {\it initial function}:
\begin{equation}\label{1.1-2}
y(x)=f(x), \quad x\in(-a,0);
\end{equation}

(ii) To assume that $q(x)=0$ a.e. on $(0,a).$\\
Actually, both ways are deeply connected. Indeed, rewriting equation (\ref{1.1-1}) in the form
$$
-y''(x)+q^+(x)y(x-a)=\lambda y(x)-r(x), \quad 0<x<1,
$$
where $r(x)=q^-(x)f(x-a)$ and
$$
q^-(x)=\left\{\begin{array}{cl}q(x), &x\in(0,a),\\[3mm]
0, & x\in(a,1),\end{array}\right. \quad
q^+(x)=\left\{\begin{array}{cl}0, &x\in(0,a), \\[3mm]
q(x), & x\in(a,1),\end{array}\right.
$$
shows that (i) leads to a non-homogenous equation, while (ii) deals with the corresponding homogenous one. That is why,
for posing an eigenvalue problem, it is natural to choose (ii). Meanwhile, (i) also may be appropriate for that purpose but
one should deal with a ''{\it linear initial function}'', i.e. when $f$ is linearly dependent on $y,$ e.g.:
$$
f(x)=y(0)g(x), \quad x\in(-a,0).
$$
Under such settings, however, one deals with the so-called frozen argument (see, e.g., \cite{ButHu21}):
$$
-y''(x)+q^+(x)y(x-a)+p(x)y(0)=\lambda y(x), \quad 0<x<1, \quad p(x):=q^-(x)g(x-a),
$$
or with other more complicated equations that are beyond the goals of our present study.

Returning to equation (\ref{1.1-1}), we note that, as in the local case $a=0,$ in order to determine a
concrete solution, one should specify also some initial conditions
$y(0)=y_0$ and $y'(0)=y_1,$
which are not supposed to be connected anyhow with the initial function in (\ref{1.1-2}). Of course, assuming additionally that
$f(x)\in C^j[-a,0]$ and $f^{(\nu)}(0)=y_\nu,\;\nu=\overline{0,j},$ for some $j\in\{0,1\},$ one would arrive at
the inclusion $y(x)\in C^j[-a,1],$
but it is not always required and may be even {\it ineligible}. The latter will be especially the case when we proceed with graphs.

\medskip
To the best of our knowledge, the first attempt of defining operator (\ref{1.1}) on a graph was made in \cite{WangYang-21, WangYang-22},
where a star-type graph was considered. Instead of (\ref{1.1-0}), the corresponding boundary value problem involved the
equations
\begin{equation}\label{1.2}
-y_j''(x)+q_j(x)y_j(x-a_j)=\lambda y_j(x), \quad 0<x<1, \quad j=\overline{1,m},
\end{equation}
but with the same matching conditions (\ref{1.3}) and boundary conditions (\ref{1.4}).  For $j=\overline{1,m},$ the functions
$q_j(x)$ belonged to $L_2(0,1),$ and $q_j(x)=0$ a.e. on $(0,a_j),$ while $a_j\in[0,1].$

Thus, relations (\ref{1.3}), (\ref{1.4}) and (\ref{1.2}) can be called nonlocal quantum graph if not all $a_j$ vanish. Such settings,
however, could be classified as {\it locally} nonlocal because the delay on each edge does not affect the other edges. Mathematically, this
means that each equation in (\ref{1.2}) can be considered separately from all other equations.

\medskip
In the present paper, we suggest another concept of operators with delay on graphs that can be characterized as {\it globally}
nonlocal, when the delay extends through vertices of the graph.
This turns out to be inherent to operator (\ref{1.1}) even on an interval. In order to illustrate this,
we consider the boundary value problem
\begin{equation}\label{1.5}
-y''(x)+q(x)y(x-a)=\lambda y(x), \quad 0<x<2, \quad y(0)=y(2)=0,
\end{equation}
where $a\in(0,2)$ is fixed and $q(x)=0$ on $(0,a).$
We note that the problem (\ref{1.5}) can be interpreted as a globally nonlocal quantum graph
with two edges. Indeed, having put
$$
y_1(x):=y(x), \quad y_2(x):=y(x+1), \quad q_1(x):=q(x), \quad q_2(x):=q(x+1), \quad 0<x<1,
$$
we consider the relations
\begin{equation}\label{1.6}
-y_j''(x)+q_j(x)y_j(x-a)=\lambda y_j(x), \quad 0<x<1, \quad j=1,2,
\end{equation}
\begin{equation}\label{1.7}
y_1(1)=y_2(0), \quad y_1'(1)=y_2'(0), \quad y_1(0)=y_2(1)=0,
\end{equation}
where $y_j(x-a)$ partially requires an additional definition. Specifically, since $q_1(x)=0$ on the interval $(0,\min\{a,1\})$
and $q_2(x)=0$ on $(0,\max\{0,a-1\}),$
the function $y_j$ in (\ref{1.6}) for negative arguments is actually absent when $j=1,$ while for $j=2,$ we define it by
\begin{equation}\label{1.7-1}
y_2(x-a):=y_1(x-a+1), \quad \max\{0,a-1\}<x<\min\{a,1\}.
\end{equation}
Note that (\ref{1.7-1}) specifies the initial function for the second equation in (\ref{1.6}) and makes the system
of relations (\ref{1.6}) and (\ref{1.7}) equivalent to the problem (\ref{1.5}). On the other hand, it can be
viewed as an eigenvalue problem on a 2-star graph illustrated on Fig.~2.
\begin{center}
\unitlength=1.0mm
\begin{picture}(100,17)
 \put(0,10){\line(1,0){100}}

 \put(0,10){\circle*{1.5}}
 \put(50,10){\circle*{1.5}}
 \put(100,10){\circle*{1.5}}

 \put (23,6.5){\small $e_1$}
 \put (74,6.5){$e_2$}

 \put (-2,5){\small $v_0$}
 \put (48.2,5){\small $v_1$}
 \put (98.3,5){\small $v_2$}

 \put (0,13){\small $0$}
 \put (46,13){\small $1$}
 \put (50.5,13){\small $0$}
 \put (96,13){\small $1$}

 \put (35,-2){\small Fig. 2. Graph $\Gamma_2$}
\end{picture}
\end{center}
The fact that equation (\ref{1.6}) for $j=2$ (i.e. on the edge $e_2)$ always involves the unknown function $y_1$
on $e_1$ means that the delay "passes" through the internal vertex $v_1.$

The illustrated idea can be directly generalized to an $m$-star graph:
\begin{equation}\label{1.8}
-y_j''(x)+q_j(x)y_j(x-a)=\lambda y_j(x), \quad 0<x<1, \quad j=\overline{1,m},
\end{equation}
\begin{equation}\label{1.9}
y_j(x-a)=y_1(x-a+1), \quad \max\{0,a-1\}<x<\min\{a,1\}, \quad j=\overline{2,m},
\end{equation}
\begin{equation}\label{1.10}
y_1(1)=y_2(0)=\ldots= y_m(0), \quad y'_1(1)=\sum_{j=2}^m y_j'(0),
\end{equation}
\begin{equation}\label{1.11}
y_j^{(\nu_j)}(1-\delta_{j,1})=0, \quad \nu_j\in\{0,1\}, \quad j=\overline{1,m},
\end{equation}
where $\delta_{1,j}$ is the Kronecker delta. We also assume that $q_j(x)\in L_2(0,1)$ for $j=\overline{1,m},$ while $q_1(x)=0$ a.e. on
$(0,\min\{a,1\})$ and $q_j(x)=0$ a.e. on $(0,\max\{0,a-1\})$ for $j=\overline{2,m}.$
Denote the corresponding graph by $\Gamma_m$ (see Fig.~3), which differs from the graph $G_m$ only by the names of two vertices and the
parametrization of all but one edges. The latter means, in particular, that matching conditions (\ref{1.10}) are equivalent to (\ref{1.3}).
The new-type conditions (\ref{1.9}) can be referred to as {\it initial-function conditions} because they actually specify an
initial function for all edges except the first one.
\begin{center}
\unitlength=0.7mm
\begin{picture}(80,90)
 \put(40,46){\line(2,3){27.7}}
 \put(40,46){\line(-1,0){50}}
 \put(40,46){\line(5,2){46.65}}
 \multiput(40,46)(2.5,-1){19}{\circle*{0.7}}
 \put(40,46){\line(2,-3){27.7}}

 \put(-10,46){\circle*{1}}
 \put(40,46){\circle*{1}}
 \put(68,88){\circle*{1}}
 \put(87,64.8){\circle*{1}}
 \put(68,4){\circle*{1}}

 \put (11,41){\small $e_1$}
 \put (55,67){\small $e_2$}
 \put (64,52.5){\small $e_3$}
 \put (47,21){\small $e_m$}

 \put (-18,44){\small $v_0$}

 \put (35,41){\small $v_1$}

 \put (68,90){\small $v_2$}
 \put (88,66){\small $v_3$}
 \put (68,0){\small $v_m$}

 \put (-10,48){\small $0$}
 \put (34,48){\small $1$}
 \put (38,51){\small $0$}
 \put (44,49.8){\small $0$}
 \put (42.34,42.0){\small $0$}
 \put (61,86){\small $1$}
 \put (82,66){\small $1$}
 \put (65,8){\small $1$}

 \put (-10,-3){\small Fig. 3. Graph $\Gamma_m$}

\end{picture}
\end{center}

In principle, this idea can be generalized to graphs of any possible structure. In Section~7, we introduce operators with global delay on an
arbitrary compact tree with nonequal edges. Various graphs with cycles can be covered too but in somewhat more complicated ways.

As was already mentioned, another aim of the present paper is to study an inverse spectral problem under the globally nonlocal
settings. For simplicity, we restrict ourself to the simplest non-degenerate star graph (i.e. the graph possessing 3 edges)
and take $a=1.$ For this case, however, we address a wide range of questions usually raised in the inverse spectral theory.

\medskip
The paper is organized as follows. In the next section, we remind some related facts and formulate our results for the inverse problems.
In Section~3, we introduce the so-called global general solution on the graph. In Section~4, we construct
the characteristic determinants of related eigenvalue problems, whose properties will be thoroughly investigated in Section~5.
In Section~6, we prove the inverse problem results formulated in Section~2.
In Section~7, we introduce operators with global delay on arbitrary compact trees and provide some examples.
\\

{\large\bf 2. Inverse problems}
\\

Inverse problems of spectral analysis consist in recovering operators from their spectral characteristics. The first substantial study
in this direction was carried out by Borg \cite{B}, who proved that the potential $q(x)$ in the Sturm--Liouville equation
\begin{equation}\label{0.1}
-y''(x)+q(x)y(x)=\lambda y(x), \quad 0<x<\pi,
\end{equation}
is uniquely determined by specifying the spectra of two boundary value problems for (\ref{0.1}) with one common boundary condition, e.g.:
$$
y(0)=y^{(\nu)}(\pi)=0, \quad \nu=0,1.
$$
Borg also proved local solvability of this inverse problem.
Although \cite{B} dealt with real $q(x),$ analogous results hold also in the non-self-adjoint case \cite{BK19}. The
self-adjoint case, however, admits a characterization of the spectra in terms of their location and asymptotics \cite{MO}.

Borg's statement of the inverse problem became a prototype of the future statements for more
complicated classes of operators. In particular, Yurko \cite{Yur-05} established that for the unique recovery of all Sturm--Liouville
equations on a compact tree, it is sufficient to know $n$ spectra, where $n$ is the number of boundary vertices.
Although this formally generalizes Borg's unique\-ness theorem on an interval, no solvability results for Yurko's statement of the inverse
problem on trees still exist.
Bondarenko \cite{Bond-20} suggested to consider an artificially overdetermined statement of the inverse problem on trees, which allows,
however, to obtain necessary and sufficient conditions for its solvability by reducing to the matrix Sturm--Liouville operator.

Concerning operators with delay as well as other classes of nonlocal operators, they require methods and approaches going beyond the
classical inverse spectral theory. Related results often reveal also qualitative differences from inverse problems for
purely differential operators.
For example, while the uniqueness theorem
for Borg-type inverse problems involving the equation
$$
-y''(x)+q(x)y(x-a)=\lambda y(x), \quad 0<x<\pi,
$$
holds if $a\in[2\pi/5,\pi)$ (see, e.g., \cite{ButYur19,BondYur18-1,VPV19,DV19}), the recent papers \cite{DB21, DB21-2, DB22} establish
its failure for $a\in(0,2\pi/5).$ Other details of the inverse spectral theory for Sturm--Liouville-type operators with constant delay
can be found, e.g., in a brief survey provided in \cite{BMSh21}.

In regard to inverse problems for operators with delay on a graph, the single study in this direction
\cite{WangYang-22} addressed the locally nonlocal case and dealt with the uniqueness and an algorithm for a partial inverse problem.
The latter consisted in recovering $q_1(x)$ and $a_1$ from an appropriate subspectrum of the problem (\ref{1.3}), (\ref{1.4})
and (\ref{1.2}) provided that $q_j(x)$ and $a_j$ for $j=\overline{2,m}$ were known {\it a priori}.
It was also assumed that $a_j\in[1/2,1)$ for $j=\overline{1,m}.$

Here, we study not partial but complete inverse problem under our globally nonlocal settings. For simplicity, we restrict ourself to
a graph with three edges and take $a=1.$ However, besides a uniqueness theorem (Theorem~2 below) and constructive procedures
for solving the inverse problem (Algorithms~1 and~2 in Section~6), we obtain also its global solvability (Theorems~3 and~4). Moreover, we
establish the uniform stability of the inverse problem (Theorem~5), which belongs to results of a new type even for classical inverse
problems on an interval.

\medskip
For $k=1,2,$ denote by ${\cal G}_k(q_2,q_3)$ the eigenvalue problem (\ref{1.8})--(\ref{1.11}) under the settings
$$
m=3, \quad a=1, \quad q_1=0, \quad \nu_1=0, \quad \nu_2=k-1, \quad \nu_3=2-k.
$$
In other words, the problem ${\cal G}_k(q_2,q_3)$ has the form
\begin{equation}\label{1.13}
-y_1''(x)=\lambda y_1(x), \quad 0<x<1,
\end{equation}
\begin{equation}\label{1.14}
-y_j''(x)+q_j(x)y_1(x)=\lambda y_j(x), \quad 0<x<1, \quad j=2,3,
\end{equation}
\begin{equation}\label{1.15}
y_1(1)=y_2(0)=y_3(0), \quad y_1'(1)=y_2'(0)+y_3'(0),
\end{equation}
\begin{equation}\label{1.16}
y_1(0)=y_2^{(k-1)}(1)=y_3^{(2-k)}(1)=0.
\end{equation}
Denote by $\Lambda_k$ the spectrum of the problem ${\cal G}_k(q_2,q_3)$ and put
\begin{equation}\label{1.16.0}
\omega_j=\int_0^1 q_j(x)\,dx, \quad j=2,3.
\end{equation}
For briefness, we agree that one and the same symbol $\{\varkappa_n\}$ will denote {\it different} sequences in~$l_2$ for all ranges of $n$
considered. Let $\{(z_{n,k}^0)^2\}_{n\ge0}$ be the sequence of zeros with account of multiplicity of the entire function
$F_k(\lambda)$ determined by the formula
\begin{equation}\label{1.16.2}
F_k(\lambda)=S_k(\sqrt\lambda), \quad S_k(\rho):=\rho\frac{1+3\cos2\rho}2\sin\rho
-\frac{\omega_2+\omega_3}4\cos2\rho+(-1)^k\frac{\omega_2-\omega_3}4.
\end{equation}
For definiteness, we assume that $0\le {\rm Re\,}z_{n,k}^0\le {\rm Re\,}z_{n+1,k}^0$ for $n\ge0.$

\medskip
{\bf Theorem 1. }{\it For $k=1,2,$ the spectrum $\Lambda_k$ consists of infinitely many complex eigenvalues having the form
\begin{equation}\label{1.16.1}
\lambda_{n,k}=z_{n,k}^2, \quad z_{n,k}=z_{n,k}^0+\frac{\varkappa_n}n, \quad n\ge1,
\end{equation}
i.e. $\Lambda_k=\{\lambda_{n,k}\}_{n\in{\mathbb N}}$ with account of multiplicity.}

\medskip
The proof of Theorem~1 as well as the following its corollary will be given in Section~5.

\medskip
{\bf Corollary 1. }{\it For $k=1,2,$ the spectrum $\Lambda_k$ can be subdivided into two subspectra:
\begin{equation}\label{1.17}
\Lambda_k=\{\mu_{n,k}\}_{|n|\in{\mathbb N}}\cup\{\xi_{n,k}\}_{n\ge0}
\end{equation}
with account of multiplicity that have the forms
\begin{equation}\label{1.19}
\mu_{n,k}=\eta_{n,k}^2, \quad \eta_{n,k}=\pi n +\sigma +(-1)^n\frac{\gamma_k}{\pi n} + \frac{\varkappa_n}n,
\end{equation}
\begin{equation}\label{1.18}
\xi_{n,k}=\zeta_{n,k}^2, \quad \zeta_{n,k}=\pi n +(-1)^n\frac{\omega_{k+1}}{4\pi n} +\frac{\varkappa_n}n \quad (when\;\;n\ne0),
\end{equation}
where
\begin{equation}\label{1.18.1}
\sigma=\frac12\arccos\Big(-\frac13\Big), \quad \gamma_k= \frac{\sqrt3}{16}\Big(\frac{\omega_2+\omega_3}3+(-1)^k(\omega_2-\omega_3)\Big).
\end{equation}
}

Although Corollary~1 does not forbid the subspectra $\{\mu_{n,k}\}_{|n|\in{\mathbb Z}}$ and $\{\xi_{n,k}\}_{n\ge0}$ to intersect, each
spectrum $\Lambda_k,$ due to (\ref{1.19}) and (\ref{1.18}), has at most finite number of multiple eigenvalues.

Consider the following inverse problem.

\medskip
{\bf Inverse Problem 1.} Given $\Lambda_1$ and $\Lambda_2,$ find $q_2(x)$ and $q_3(x).$

\medskip
As in the case of an interval (see \cite{ButYur19, ButDjur22}), for the unique determination of the potentials, it is sufficient
to specify appropriate subspectra. For example, the following theorem holds.

\medskip
{\bf Theorem 2. }{\it Specification of $\{\mu_{n,k}\}_{|n|\in{\mathbb N}},\;k=1,2,$ uniquely determines $q_2(x)$ and $q_3(x).$}

\medskip
This theorem as well as Theorems 3--5 formulated below will be proved in Section~6.

In what follows, we refer to $\{\mu_{n,k}\}_{|n|\in{\mathbb N}}$ as $\mu$-{\it subspectrum}. Thus, Theorem~2 is a uniqueness theorem for
the following inverse problem.

\medskip
{\bf Inverse Problem 2.} Given $\mu$-subspectra $\{\mu_{n,k}\}_{|n|\in{\mathbb N}},\,k=1,2,$ find $q_2(x)$ and $q_3(x).$

\medskip
Unlike Inverse Problem 1, it is not overdetermined, which follows from the next theorem.

\medskip
{\bf Theorem 3. }{\it For arbitrary complex sequence $\{\mu_{n,1}\}_{|n|\in{\mathbb N}}$ and $\{\mu_{n,2}\}_{|n|\in{\mathbb N}}$ of the form
(\ref{1.19}) with any complex numbers $\gamma_1$ and $\gamma_2,$ there exist functions $q_2(x)$ and $q_3(x)$ in $L_2(0,1)$
such that these sequences are $\mu$-subspectra of the problems ${\cal G}_1(q_2,q_3)$ and ${\cal G}_2(q_2,q_3).$}

\medskip
Thus, the asymptotics in (\ref{1.19}) form necessary and sufficient conditions for the solvability of Inverse Problem~2. In other words,
(\ref{1.19}) gives a complete characterization of $\mu$-subsectra.

Let us return to Inverse Problem~1. In spite of the non-minimality of its input data, one can formulate necessary and sufficient conditions
for its solvability analogously to how it was made for a pencil with two delays \cite{BMSh21} or for a Dirac-type system with delay
\cite{ButDjur22}.

\medskip
{\bf Theorem 4. }{\it For arbitrary complex sequences $\{\lambda_{n,1}\}_{n\in{\mathbb N}}$
and $\{\lambda_{n,2}\}_{n\in{\mathbb N}}$ to be the spectra of the boundary
value problems ${\cal G}_1(q_2,q_3)$ and ${\cal G}_2(q_2,q_3)$ with some functions $q_2(x)$ and $q_3(x)$ in $L_2(0,1),$ it is
necessary and sufficient to satisfy the following two conditions:

(i) For $k=1,2,$ the sequence $\{\lambda_{n,k}\}_{n\in{\mathbb N}}$ has the form~(\ref{1.16.1}), while $\{(z_{n,k}^0)^2\}_{n\ge0}$ are zeros
of the function $F_k(\lambda)$ determined by (\ref{1.16.2}) with arbitrary fixed $\omega_2,\omega_3\in{\mathbb C};$

(ii) The exponential types of the entire functions $(\Delta_k-\Delta_0)(\rho^2),$ $k=1,2,$ in the $\rho$-plane do not exceed~$2,$ where the
functions $\Delta_\nu(\lambda)$ are determined by the formulae
\begin{equation}\label{1.19.1}
\Delta_0(\lambda):=\frac{\sin\sqrt\lambda}{2\sqrt\lambda}(1+3\cos2\sqrt\lambda),
\quad \Delta_k(\lambda):=-\frac{\alpha_k}{\lambda_{0,k}^0}\prod_{n\in{\mathbb
N}}\frac{\lambda_{n,k}-\lambda}{\lambda_{n,k}^0}, \quad k=1,2,
\end{equation}
where, in turn,
$$
\alpha_k=\lim_{\lambda\to0}\frac{F_k(\lambda)}{\lambda^{s_k}}, \quad
\lambda_{n,k}^0=\left\{\begin{array}{cl} (z_{n,k}^0)^2, & z_{n,k}^0\ne0, \quad \\[3mm]
-1, & z_{n,k}^0=0, \end{array}\right. \quad n\ge0,
$$
while $s_k$ is the multiplicity of a zero of $F_k(\lambda)$ at the origin. }

\medskip
Finally, let us formulate the uniform stability for Inverse Problem~1. For this purpose, along with the problems ${\cal G}_k(q_2,q_3),$ we
will consider problems ${\cal G}_k(\tilde q_2,\tilde q_3)$ of the same forms but with other potentials $\tilde q_2(x), \tilde
q_3(x)\in L_2(0,1).$ If some symbol $\gamma$ denotes an object related to ${\cal G}_k(q_2,q_3),$ then this symbol with tilde $\tilde\gamma$
will denote the analogous object related to ${\cal G}_k(\tilde q_2,\tilde q_3).$

\medskip
{\bf Theorem 5. }{\it For any $r>0,$ there exists $C_r>0$ such that the estimate
\begin{equation}\label{1.20}
\sum_{j=2}^3\|q_j-\tilde q_j\|_{L_2(0,1)}\le C_r\sum_{k=1}^2\|\{n(z_{n,k}-\tilde z_{n,k})\}_{n\in{\mathbb N}}\|_{l_2}
\end{equation}
holds whenever $\|\{n(z_{n,k}-z_{n,k}^0)\}_{n\in{\mathbb N}}\|_{l_2}\le r$
and $\|\{n(\tilde z_{n,k}-z_{n,k}^0)\}_{n\in{\mathbb N}}\|_{l_2}\le r$ for $k=1,2.$ }

\medskip
Although stability results in the inverse spectral theory go back to Borg \cite{B}, the uniform stability of the classical inverse
Sturm--Liouville problem for the first time was established in \cite{SavShk}. In \cite{But21-1, But21-2},
this type of stability was addressed for some classes of integro-differential operators, for which a different from \cite{SavShk}
approach was used, as an integral part of which there became proving the uniform stability of recovering the characteristic function of an
operator under consideration from its zeros, i.e. from the spectrum of this operator. In \cite{But22}, the latter type results were
extended to a more general class of entire functions. In \cite{ButDjur22}, those results were used for proving the uniform stability
of an inverse problem for Dirac-type operators with constant delay. Here, we use them for proving Theorem~5. We also obtain other
properties of the characteristic functions as corollaries from the corresponding general assertions in \cite{But22}.

Theorem~5 means that Inverse Problem~1 is Lipschitz continuous on each ball of a finite radius. As in \cite{But21-1,
But21-2}, one can obtain an analogous estimate also in the $\lambda$-plane. Moreover, since, according to Theorem~2, it is sufficient to
specify only $\mu$-subspectra for the unique recovery of $q_2(x)$ and $q_3(x),$ one can expect that some analogue of Theorem~5 will hold
also for Inverse Problem~2. However, obtaining precise statements is beyond the present study.

Results analogous to Theorems~1--5 can be obtained also in the case $a\ne1,$ appearing more difficult for $a<1,$ when the dependence of the
characteristic function on the potentials becomes nonlinear. The nonlinear case may also require an additional spectrum to be specified.
Moreover, one can expect the failure of the corresponding uniqueness theorem for small positive~$a$ as in the case of an interval
\cite{DB21, DB21-2, DB22}.
Analogously, one can study operators with global delay on non-star graphs (see Section~7) starting with illustrative Examples~2--4.
\\

{\large\bf 3. Global general solution}
\\

The investigation of purely differential operators on graphs usually begins with constructing appropriate fundamental systems of solutions
for the corresponding differential equation on each edge. However, it is not possible in our globally nonlocal settings because any
solution on each edge (except the first one) depends on a solution (or solutions) on another edge (edges).
This makes it necessary to construct a {\it global} general solution (GGS) on the entire graph.

Any function $y$ on a graph (not necessarily star graph) can be interpreted as the vector
\begin{equation}\label{3.1}
y=[y_1,y_2,\ldots,y_m],
\end{equation}
where $y_j(x)$ is defined on the $j$-th edge, while $m$ is the number of edges in this graph. Since the edge lengths may differ,
(\ref{3.1}) is, in general, a vector of functions rather than vector-function. In particular, for the star graph $\Gamma_m$ on Fig.~3,
the function $y_j(x)$ is defined on the edge~$e_j.$

By the GGS on $\Gamma_m,$ we mean a function $y$ in sense of (\ref{3.1}) whose $j$-th component $y_j(x)$ solves the $j$-th
equation in (\ref{1.8}) and obeys the $j$-th initial-function condition (\ref{1.9}) if $j\ne1,$ while neither the matching conditions
(\ref{1.10}) nor the boundary conditions (\ref{1.11}) are supposed to be fulfilled. As will be seen below, such solution $y$ involves $2m$
indefinite constants.

Here, we aim to construct the GGS on $\Gamma_m.$ Analogously, GGS can be constructed also
on an arbitrary globally nonlocal tree introduced in Section~7.

First, let $a\in(0,1].$ For $j=\overline{1,m},$ denote
$$
q_j^-(x):=\left\{\begin{array}{cl}q_j(x), &x\in(0,a),\\[3mm]
0, & x\in(a,1),\end{array}\right. \quad
q_j^+(x):=\left\{\begin{array}{cl}0, &x\in(0,a),\\[3mm]
q_j(x), & x\in(a,1).\end{array}\right.
$$
Then the $j$-th equation in (\ref{1.8}) takes the form
\begin{equation}\label{3.3}
-y_j''(x)+q_j^+(x)y_j(x-a)=\lambda y_j(x)-f_j(x), \quad 0<x<1,
\end{equation}
where, with accordance to (\ref{1.9}), we have
\begin{equation}\label{3.4}
f_j(x)=q_j^-(x)y_1(x-a+1).
\end{equation}
Consider the fundamental system of solutions $\{Y_{1,j}(x,\lambda),Y_{2,j}(x,\lambda)\}$
of the homogenous equation corresponding to (\ref{3.3}) (i.e. without~$f_j)$
under the initial conditions $Y_{\nu,j}^{(l)}(0,\lambda)=\delta_{\nu,l+1},$ $l=0,1.$ These solutions can be constructed by the explicit
formula (see, e.g., \cite{DB22})
\begin{equation}\label{3.5}
Y_{\nu,j}(x,\lambda)=\sum_{k=0}^N Y_{\nu,j,k}(x,\lambda), \;\;
Y_{\nu,j,k}(x,\lambda)=\int\limits_{ka}^x\frac{\sin\rho(x-t)}\rho q_j^+(t)Y_{\nu,j,k-1}(t-a,\lambda)\,dt, \;\; k\ge1,
\end{equation}
where $N\in{\mathbb N}$ is such that $a\in[1/(N+1),1/N)$ and  $\rho^2=\lambda,$ while
\begin{equation}\label{3.6}
Y_{1,j,0}(x,\lambda)=\cos\rho x, \quad Y_{2,j,0}(x,\lambda)=\frac{\sin\rho x}\rho.
\end{equation}
Note that the function
\begin{equation}\label{3.7}
y_1(x)=C_{1,1}Y_{1,1}(x,\lambda) +C_{2,1}Y_{2,1}(x,\lambda)
\end{equation}
is the general solution on the edge $e_1,$ being, in turn, a unique edge on which any solution is independent of solutions on any other edges.

Further, for $j=\overline{2,m},$ denote by $z_j(x,\lambda)$ the unique solution of the non-homogenous equation (\ref{3.3})
with the free term $f_j(x)$ determined by (\ref{3.4}) and
under the homogenous initial conditions $z_j(0,\lambda)=z_j'(0,\lambda)=0.$ By virtue of (\ref{3.4}) and (\ref{3.7}), we have
$$
z_j(x,\lambda)=C_{1,1}Z_{1,j}(x,\lambda) +C_{2,1}Z_{2,j}(x,\lambda),
$$
where, the function $Z_{\nu,j}(x,\lambda)$ for $\nu=1,2$ is the unique solution of the equation (\ref{3.3}) with the free term
$f_j(x)=q_j^-(x)Y_{\nu,1}(x-a+1,\lambda)$ under the conditions $Z_{\nu,j}(0,\lambda)=Z_{\nu,j}'(0,\lambda)=0.$

Analogously to (\ref{3.5}), one can obtain the following explicit formula for $\nu=1,2:$
\begin{equation}\label{3.8}
Z_{\nu,j}(x,\lambda)=\sum_{k=0}^N Z_{\nu,j,k}(x,\lambda), \;\;
Z_{\nu,j,k}(x,\lambda)=\int\limits_{ka}^x\frac{\sin\rho(x-t)}\rho q_j^+(t)Z_{\nu,j,k-1}(t-a,\lambda)\,dt, \;\; k\ge1,
\end{equation}
with the same $N,$ where
\begin{equation}\label{3.9}
Z_{\nu,j,0}(x,\lambda)=\int_0^x\frac{\sin\rho(x-t)}\rho q_j^-(t)Y_{\nu,1}(t-a+1,\lambda)\,dt.
\end{equation}

Thus, the required GGS in the case $a\in(0,1]$ has the form (\ref{3.1}), where $y_1$ is determined by (\ref{3.7}),
while the rest components have the form
$$
y_j(x)=C_{1,j}Y_{1,j}(x,\lambda) +C_{2,j}Y_{2,j}(x,\lambda) +C_{1,1}Z_{1,j}(x,\lambda) +C_{2,1}Z_{2,j}(x,\lambda), \quad
j=\overline{2,m}.
$$
Here, $C_{\nu,j}$ for $\nu=1,2$ and $j=\overline{1,m}$ are arbitrary constants, while the functions $Y_{\nu,j}(x,\lambda)$ and
$Z_{\nu,j}(x,\lambda)$ are determined by the explicit formulae (\ref{3.5}), (\ref{3.6}) and  (\ref{3.8}), (\ref{3.9}), respectively.

For $a\in[1,2),$ the first component of the GGS will, obviously, take the form
\begin{equation}\label{3.11}
y_1(x)=C_{1,1}\cos\rho x +C_{2,1}\frac{\sin\rho x}\rho.
\end{equation}
According to (\ref{1.9}) along with the requirement that $q_j(x)=0$ a.e. on $(0,a-1)$ as $j=\overline{2,m},$ the rest equations in
(\ref{1.8}) can be rewritten as
$$
-y_j''(x)+q_j(x)y_1(x-a+1)=\lambda y_j(x), \quad 0<x<1, \quad j=\overline{2,m}.
$$
Hence, by virtue of (\ref{3.11}), we have
\begin{equation}\label{3.12}
y_j(x)=C_{1,j}\cos\rho x +C_{2,j}\frac{\sin\rho x}\rho +C_{1,1}P_j(x,\lambda) +C_{2,1}Q_j(x,\lambda), \quad
j=\overline{2,m},
\end{equation}
where
\begin{equation}\label{3.13}
P_j(x,\lambda)=\int_{a-1}^x\frac{\sin\rho(x-t)}\rho q_j(t)\cos\rho(t-a+1)\,dt,
\end{equation}
\begin{equation}\label{3.14}
Q_j(x,\lambda)=\int_{a-1}^x \frac{\sin\rho(x-t)}\rho q_j(t)\frac{\sin\rho(t-a+1)}\rho\,dt.
\end{equation}
Thus, the GGS in the case $a\in[1,2)$ is determined by (\ref{3.1}) along with the formulae (\ref{3.11})--(\ref{3.14}).

The described scheme of constructing the GGS can be extended in the obvious way to an arbitrary tree (see Section~7).

We note finally that substituting the obtained GGS into the $m$ matching conditions (\ref{1.10}) and the $m$ boundary conditions
(\ref{1.11}) gives a homogenous algebraic system with respect to the $2m$ undetermined constants $C_{\nu,j}.$ Eigenvalues
of the corresponding boundary value problem will coincide with zeros of its determinant, which is called
{\it characteristic determinant} or {\it characteristic function}. In the next section, we construct
the characteristic functions of the boundary value problems introduced in the preceding section.
\\

{\large\bf 4. Characteristic functions of the problems ${\cal G}_k(q_2,q_3)$}
\\

It is easy to check that the components
\begin{equation}\label{2.1}
y_1(x)=C_1\frac{\sin\rho x}\rho, \quad y_j(x)=C_{1,j}\cos\rho x +C_{2,j}\frac{\sin\rho x}\rho +C_1Q_j(x,\lambda), \quad j=2,3,
\end{equation}
where
\begin{equation}\label{2.2}
Q_j(x,\lambda)=\int_0^x \frac{\sin\rho(x-t)}\rho q_j(t)\frac{\sin\rho t}\rho\,dt,
\end{equation}
form the GGS $y=[y_1,y_2,y_3]$  of the system that consists of (\ref{1.13})
and (\ref{1.14}) under the standing Dirichlet boundary condition  $y_1(0)=0.$
Hence, for $y$ to be an eigenfunction of the problem ${\cal G}_k(q_2,q_3),$
besides being nonzero, it remains to satisfy the rest boundary conditions in (\ref{1.16}) as well as all matching
conditions~(\ref{1.15}).
Thus, substituting (\ref{2.1}) and (\ref{2.2}) into (\ref{1.15}) and (\ref{1.16}), we conclude that
$y$ is an eigenfunction if and only if the column-vector ${\cal
C}:=[C_1,C_{1,2},C_{2,2},C_{1,3},C_{2,3}]^T$ is nonzero and it satisfies the linear algebraic system ${\cal A}_k(\lambda){\cal C}=0,$ where,
in particular,
$$
{\cal A}_1(\lambda)=\left(\begin{array}{ccccc} \displaystyle
\frac{\sin\rho}\rho  &     -1     &                    0                  &         0      &      0      \\[5mm]
          0          &     -1     &                    0                  &         1      &      0      \\[5mm]
      \cos\rho       &      0     &                   -1                  &         0      &     -1      \\[5mm]
  Q_2(1,\lambda)     &  \cos\rho  &  \displaystyle \frac{\sin\rho}\rho    &         0      &      0      \\[5mm]
  Q_3'(1,\lambda)    &      0     &                    0                  &  -\rho\sin\rho &   \cos\rho
\end{array}\right).
$$
Here and below, $f'$ denotes the derivative of any function $f$ with respect to its {\it first} argument. Thus, eigenvalues of the problem
${\cal G}_1(q_2,q_3)$ coincide with zeros of the determinant
$$
\Delta_1(\lambda):=\det{\cal A}(\lambda)=\Delta_0(\lambda) +Q_2(1,\lambda)\cos\rho +Q_3'(1,\lambda)\frac{\sin\rho}\rho,
$$
where $\Delta_0(\lambda)$ is defined in (\ref{1.19.1}). The function $\Delta_1(\lambda)$ is called {\it characteristic determinant} or {\it
characteristic function} of the problem ${\cal G}_1(q_2,q_3).$

Symmetrically, the characteristic function of ${\cal G}_k(q_2,q_3)$ for $k=1,2$ has the form
\begin{equation}\label{2.5}
\Delta_k(\lambda)=\Delta_0(\lambda) +Q_{k+1}(1,\lambda)\cos\rho +Q_{4-k}'(1,\lambda)\frac{\sin\rho}\rho.
\end{equation}
For $j=2,3,$ formula (\ref{2.2}) implies the representations
\begin{equation*}\label{2.6}
Q_j(x,\lambda)=\frac1{2\rho^2}\int_0^x q_j(t)\cos\rho(x-2t)\,dt -\frac{\omega_j(x)}{2\rho^2}\cos\rho x, \quad \omega_j(x):=\int_0^x
q_j(t)\,dt,
\end{equation*}
\begin{equation*}\label{2.7}
Q_j'(x,\lambda)=\frac{\omega_j(x)}{2\rho}\sin\rho x -\frac1{2\rho}\int_0^x q_j(t)\sin\rho(x-2t)\,dt.
\end{equation*}
Substituting $x=1$ and changing the variable of integration, we arrive at
\begin{equation}\label{2.8}
Q_j(1,\lambda)=\int_0^1 u_j^+(x)\frac{\cos\rho x}{\rho^2}\,dx -\omega_j\frac{\cos\rho}{2\rho^2}, \quad
Q_j'(1,\lambda)=\omega_j\frac{\sin\rho}{2\rho} +\int_0^1 u_j^-(x)\frac{\sin\rho x}\rho\,dx,
\end{equation}
where $\omega_j=\omega_j(1)$ and
\begin{equation}\label{2.9}
u_j^+(x)=\frac14\Big(q_j\Big(\frac{1+x}2\Big)+q_j\Big(\frac{1-x}2\Big)\Big), \quad
u_j^-(x)=\frac14\Big(q_j\Big(\frac{1+x}2\Big)-q_j\Big(\frac{1-x}2\Big)\Big).
\end{equation}
Using (\ref{2.8}), we calculate
\begin{equation}\label{2.10}
Q_j(1,\lambda)\cos\rho=\int_0^1 u_j^+(1-x)\frac{\cos\rho x}{2\rho^2}\,dx + \int_1^2 u_j^+(x-1)\frac{\cos\rho x}{2\rho^2}\,dx
-\omega_j\frac{1+\cos2\rho}{4\rho^2},
\end{equation}
\begin{equation}\label{2.11}
Q_j'(1,\lambda)\frac{\sin\rho}\rho =\int_0^1 u_j^-(1-x)\frac{\cos\rho x}{2\rho^2}\,dx - \int_1^2 u_j^-(x-1)\frac{\cos\rho x}{2\rho^2}\,dx
+\omega_j\frac{1-\cos2\rho}{4\rho^2}.
\end{equation}
Substituting (\ref{2.10}) and (\ref{2.11}) into (\ref{2.5}), we arrive at the following assertion.

\medskip
{\bf Lemma 1. }{\it For $k=1,2,$ the characteristic function of ${\cal G}_k(q_2,q_3)$ has the form
\begin{equation}\label{2.12}
\Delta_k(\lambda)=\Delta_0(\lambda)+(-1)^k\frac{\omega_2-\omega_3}{4\rho^2} -\frac{\omega_2+\omega_3}{4\rho^2}\cos2\rho +\int_0^2
w_k(x)\frac{\cos\rho x}{\rho^2}\,dx,
\end{equation}
where
\begin{equation}\label{2.13}
w_k(x)=\frac12\left\{\begin{array}{l}(u_{k+1}^+ +u_{4-k}^-)(1-x), \quad x\in(0,1),\\[3mm]
(u_{k+1}^+ -u_{4-k}^-)(x-1), \quad x\in(1,2),
\end{array}\right.
\end{equation}
while $\omega_j$ and $u_j^\pm(x)$ are defined in (\ref{1.16.0}) and (\ref{2.9}), respectively.}

\medskip
Thus, for any functions $q_2(x)$ and $q_3(x)$ in $L_2(0,1),$ relations (\ref{2.9}) and (\ref{2.13}) uniquely determine a pair of functions
$w_1(x)$ and $w_2(x)$ in $L_2(0,2).$ One can easily notice that this correspondence is one-to-one. Moreover, the following assertion holds.

\medskip
{\bf Lemma 2. }{\it Specification of the functions $w_1(x),\,w_2(x)\in L_2(0,2)$ uniquely determines the potentials $q_2(x)$ and $q_3(x)$ by
the formulae
\begin{equation}\label{2.14}
\left.\begin{array}{c} q_2(x)=2(w_1-w_2)(2x)+2(w_1+w_2)(2-2x),\\[3mm]
q_3(x)=2(w_2-w_1)(2x)+2(w_1+w_2)(2-2x),
\end{array}\right\} \quad x\in(0,1).
\end{equation}
}

{\it Proof.} Rewrite relations (\ref{2.13}) in the form
$$
2w_k(1-x)=(u_{k+1}^++u_{4-k}^-)(x), \quad 2w_k(1+x)=(u_{k+1}^+-u_{4-k}^-)(x), \quad x\in(0,1), \quad k=1,2,
$$
which immediately implies
\begin{equation}\label{2.15}
\left.\begin{array}{c}
\displaystyle u_j^+(x)=w_{j-1}(1-x)+w_{j-1}(1+x),\\[3mm]
\displaystyle u_j^-(x)=w_{4-j}(1-x)-w_{4-j}(1+x),
\end{array}\right\} \quad x\in(0,1), \quad j=2,3.
\end{equation}
Further, inverting relations (\ref{2.9}), we obtain
$$
q_j\Big(\frac{1\pm x}2\Big)=2(u_j^+\pm u_j^-)(x), \quad x\in(0,1),
$$
or
\begin{equation}\label{2.16}
q_j(x)=2\left\{\begin{array}{c}
\displaystyle (u_j^+-u_j^-)(1-2x), \quad x\in\Big(0,\frac12\Big),\\[5mm]
\displaystyle (u_j^++u_j^-)(2x-1), \quad x\in\Big(\frac12,1\Big).
\end{array}\right.
\end{equation}
Finally, substituting (\ref{2.15}) into (\ref{2.16}), we arrive at (\ref{2.14}). $\hfill\Box$
\\

{\large\bf 5. A model entire function}
\\

Consider the entire function
\begin{equation}\label{4.1}
\Delta(\lambda)=\Delta_0(\lambda)+\frac{g+h\cos2\rho}{\rho^2} +\int_0^3 w(x)\frac{\cos\rho x}{\rho^2}\,dx, \quad w(x)\in L_2(0,3),
\end{equation}
which differs from (\ref{2.12}) by the widest possible support of the function under the integral. In the next section, we actually prove
that it is an only difference, i.e. any function of the form (\ref{4.1}) is the characteristic function of some problem ${\cal
G}_k(q_2,q_3)$ if
\begin{equation}\label{4.1-1}
g+h +\int_0^3 w(x)\,dx=0
\end{equation}
and $w(x)=0$ a.e. on $(2,3).$ Condition (\ref{4.1-1}) alone is necessary and sufficient for $\Delta(\lambda)$ to be entire. The analogous
condition, obviously, holds also for the characteristic function (\ref{2.12}).

In this section, we study the characteristic functions $\Delta_k(\lambda)$ via their common form (\ref{4.1}). For convenience, we obtain all
required assertions as corollaries from the corresponding facts in \cite{But22}. For this purpose, we consider the function
$\theta(\rho):=\rho^2\Delta(\rho^2),$ which has the form
\begin{equation}\label{4.2}
\theta(\rho)=S(\rho)+\int_{-3}^3 v(x)\exp(i\rho x)\,dx, \quad v(x)\in L_2(-3,3),
\end{equation}
where we have defined:
\begin{equation}\label{4.3}
S(\rho):=\rho^2\Delta_0(\rho^2)+h\cos2\rho+g, \quad v(x):=\frac12\left\{\begin{array}{cl} w(-x), & x\in(-3,0),\\[3mm]
w(x), & x\in(0,3). \end{array}\right.
\end{equation}
Thus, the function $v(x)$ is even, which will be actually used starting from Lemma~3 below.

It is easy to see that the zeros $\{\rho_n^0\}_{|n|\in{\mathbb N}}$ of the function $\Delta_0(\rho^2)$ in the $\rho$-plane can be split into
three non-intersecting sequences:
\begin{equation}\label{4.4}
\{\rho_n^0\}_{|n|\in{\mathbb N}}=
\{\rho_{n,1}^0\}_{|n|\in{\mathbb N}}\cup\{\rho_{n,2}^0\}_{n\in{\mathbb Z}}\cup\{\rho_{n,3}^0\}_{n\in{\mathbb Z}},
\quad \rho_{n,1}^0=\pi n, \quad \rho_{n,2}^0=-\rho_{-n,3}^0=\pi n+\sigma,
\end{equation}
where $\sigma$ is defined in (\ref{1.18.1}). Thus, the function $S(\rho)$ admits the representation
\begin{equation}\label{4.4.1}
S(\rho)=P_1(\rho)S_0(\rho),
\end{equation}
where $P_1(\rho)$ is a polynomial of degree 1, while $S_0(\rho)$ is a sine-type function of type $3$ with asymptotically
separated zeros.

For convenience, we remind that an entire function $f(\rho)$ of exponential type is called {\it sine-type function} (of
type $b)$ if there exist positive constants $K,$ $C_1$ and $C_2$ such that
\begin{equation}\label{4.4.2}
C_1<|f(\rho)|\exp(-|{\rm Im}\rho|b)<C_2, \quad |{\rm Im}\rho|>K.
\end{equation}

Then, by virtue of Theorem~4 in \cite{But22}, zeros $\{z_n\}_{n\in A}$ of the function $\theta(z)$
have the form
\begin{equation}\label{4.5}
z_n=z_n^0+\frac{\varkappa_n}{\mu_n}, \quad \mu_n:=\left\{\begin{array}{cl} z_n^0, & z_n^0\ne0,\\[3mm]
-1, & z_n^0=0, \end{array}\right.
\end{equation}
where $\{z_n^0\}_{n\in A}$ are zeros of the function $S(\rho).$ As the index set $A,$ it is convenient now to use $A:=\{n:
n=\pm0,\pm1,\pm2,\ldots\},$ where $0$ and $-0$ are considered as different indices.

\medskip
{\bf Remark 1.} Since $S(\rho)$ is even, the function $F(\lambda):=S(\sqrt\lambda)$ is entire too. We agree that the zeros $\{z_n^0\}_{n\in
A}$ of $S(\rho)$ are indexed in such a way that $\{(z_n^0)^2\}_{n\ge0}$ are all zeros of $F(\lambda)$ with account of multiplicity. For
definiteness, we also assume that $0\le {\rm Re\,}z_n^0\le {\rm Re\,}z_{n+1}^0,\;n\ge0.$

\medskip
{\bf Lemma 3. }{\it The function $\Delta(\lambda)$ has infinitely many zeros $\{\lambda_n\}_{n\in{\mathbb N}},$ which have the form}
\begin{equation}\label{4.5.1}
\lambda_n=z_n^2, \quad z_n=z_n^0+\frac{\varkappa_n}n.
\end{equation}

{\it Proof.} Once we deal with an even function $v(x)$ in (\ref{4.2}), the function $\theta(z)$ is even too. Therefore, without loss of
generality, one can assume that
$$
z_{-n}=-z_n, \quad n\ge0; \quad z_{-0}=z_0=0.
$$
The latter follows from the entireness of $\Delta(\lambda),$ whose zeros thus coincide with $\{z_n^2\}_{n\in{\mathbb N}}.$

Further, by virtue of the representation (\ref{4.4.1}) along with the left-hand estimate in (\ref{4.4.2}) for $f(\rho)=S_0(\rho)$
and $b=3,$ and by Rouch\'e's theorem,
zeros of $S(\rho)$ have the form
\begin{equation}\label{4.6}
z_n^0=\rho_n^0+\varepsilon_n, \quad \varepsilon_n=o(1), \quad |n|\to\infty,
\end{equation}
where $\rho_{-0}^0=\rho_0^0=0.$ Hence, by virtue of (\ref{4.4}), we have $z_n^0\asymp n$ as $|n|\to\infty.$ Therefore, ''$\mu_n$'' in the
denominator in (\ref{4.5}) can be replaced with ''$n$'' for $n\in{\mathbb N},$ which finishes the proof. $\hfill\Box$

\medskip
{\bf Corollary 2. }{\it Zeros of the function $\Delta(\lambda)$ can be subdivided into two sequences:
\begin{equation}\label{4.6.1}
\{\lambda_n\}_{n\in{\mathbb N}}=\{\mu_n\}_{n\in{\mathbb Z}}\cup\{\xi_n\}_{n\in{\mathbb N}}
\end{equation}
with account of multiplicity that have the forms
\begin{equation}\label{4.6.2}
\mu_n=\eta_n^2, \quad \eta_n=\pi n +\sigma +(-1)^n\frac{\sqrt3}{4\pi n}\Big(g-\frac{h}3\Big) +\frac{\varkappa_n}n \quad (when\;\;n\ne0),
\end{equation}
\begin{equation}\label{4.6.3}
\xi_n=\zeta_n^2, \quad \zeta_n=\pi n -(-1)^n\frac{g+h}{2\pi n}+\frac{\varkappa_n}n,
\end{equation}
where $\sigma$ is defined in (\ref{1.18.1}). }

\medskip
{\it Proof.} Let us first refine asymptotics (\ref{4.6}). Putting $H(\rho):=\rho\Delta_0(\rho^2),$ we calculate
\begin{equation}\label{4.7}
H'(\rho)=\frac{1+3\cos2\rho}2 \cos\rho -3\sin\rho\sin2\rho.
\end{equation}
Lagrange's formula gives
\begin{equation}\label{4.7.1}
H(\rho)=(\rho-\rho_n^0)H'(\rho_n^0)+O((\rho-\rho_n^0)^2), \quad \rho\to\rho_n^0.
\end{equation}
Since the function $H(\rho)$ is bounded in each fixed horizontal strip, the estimate in (\ref{4.7.1}) is uniform with respect to
$n\in{\mathbb Z}.$ Indeed, denote $C_\delta=\sup|H(\rho)|$ as $|{\rm Im}\rho|<\delta.$ Then the maximum modulus principle for analytic
functions gives
$$
|H'(\rho_n^0)|\le\max_{|\rho-\rho_n^0|=\delta}\frac{|H(\rho)|}{|\rho-\rho_n^0|}\le\frac{C_\delta}\delta, \quad
\max_{|\rho-\rho_n^0|\le\delta}\frac{|H(\rho)-(\rho-\rho_n^0)H'(\rho_n^0)|}{|\rho-\rho_n^0|^2}\le\frac{2C_\delta}{\delta^2},
$$
where the right-hand side is independent of $n.$ Further, by substituting (\ref{4.6}) into (\ref{4.7.1}), we get
\begin{equation}\label{4.7.2}
H(z_n^0)=\varepsilon_nH'(\rho_n^0)+O(\varepsilon_n^2), \quad |n|\to\infty, \quad n\in{\mathbb Z},
\end{equation}
where, in turn, according to (\ref{4.4}) and (\ref{4.7}), we have
\begin{equation}\label{4.8}
H'(\rho_{n,1}^0)=2(-1)^n, \;\; |n|\in{\mathbb N}; \quad H'(\rho_{n,2}^0)=H'(\rho_{n,3}^0)=\frac{4(-1)^{n+1}}{\sqrt3},  \;\; n\in{\mathbb Z}.
\end{equation}
On the other hand, substituting (\ref{4.6}) into (\ref{4.3}) and using the designation $H(\rho),$ we obtain
$$
z_n^0H(z_n^0)=-h\cos2z_n^0-g, \quad n\in A,
$$
which (for sufficiently large $|n|)$ along with (\ref{4.7.2}) and (\ref{4.8}) leads to the representation
$$
\varepsilon_n=-\frac{h\cos2z_n^0+g}{z_n^0(H'(\rho_n^0)+O(\varepsilon_n))} =O\Big(\frac1n\Big), \quad |n|\to\infty.
$$
Thus, using (\ref{4.6}), we refine
\begin{equation}\label{4.9}
\varepsilon_n=\varepsilon(\rho_n^0) +O\Big(\frac1{n^2}\Big), \quad |n|\to\infty, \quad \varepsilon(\rho):=-\frac{h\cos2\rho+g}{\rho
H'(\rho)}.
\end{equation}
Moreover, by virtue of (\ref{4.4}) and (\ref{4.8}), we have
\begin{equation}\label{4.10}
\varepsilon(\rho_{n,1}^0)=(-1)^{n+1}\frac{g+h}{2\pi n}, \quad \varepsilon(\rho_{n,2}^0)=\varepsilon(\rho_{n,3}^0)=(-1)^n\frac{\sqrt3}{4\pi
n}\Big(g-\frac{h}3\Big) +O\Big(\frac1{n^2}\Big), \quad |n|\to\infty.
\end{equation}
By virtue of (\ref{4.4}), (\ref{4.6}), (\ref{4.9}) and (\ref{4.10}), zeros of $S(\rho)$ can be split into three sequences:
\begin{equation}\label{4.11}
\{z_n^0\}_{n\in A}=\{z_{n,1}^0\}_{n\in A}\cup\{z_{n,2}^0\}_{n\in{\mathbb Z}}\cup\{z_{n,3}^0\}_{n\in{\mathbb Z}}
\end{equation}
of the form
\begin{equation}\label{4.12}
z_{n,1}^0=\pi n-(-1)^n\frac{g+h}{2\pi n}+O\Big(\frac1{n^2}\Big), \quad n\in A,
\end{equation}
\begin{equation}\label{4.13}
z_{n,j}^0=\pi n+(-1)^j\sigma+(-1)^n\frac{\sqrt3}{4\pi n}\Big(g-\frac{h}3\Big) +O\Big(\frac1{n^2}\Big), \quad n\in{\mathbb Z}, \quad j=2,3,
\end{equation}
as $|n|\to\infty.$ Since $S(\rho)$ is even, one can also assume without loss of generality that
\begin{equation}\label{4.14}
z_{n,1}^0=-z_{-n,1}^0, \quad n\in A, \qquad z_{n,2}^0=-z_{-n,3}^0, \quad n\in{\mathbb Z}.
\end{equation}
Thus, according to (\ref{4.5.1}) and (\ref{4.11})--(\ref{4.14}), zeros of $\Delta(\lambda)$ in the $\lambda$-plane can be split into two
sequences: $\{(z_n)^2\}_{n\in{\mathbb N}}=\{(\eta_n)^2\}_{n\in{\mathbb Z}}\cup\{(\zeta_n)^2\}_{n\in{\mathbb N}},$ where
$$
\eta_n=z_{n,2}^0+\frac{\varkappa_n}{n+1/2}, \quad \zeta_n=z_{n,1}^0+\frac{\varkappa_n}n,
$$
which along with (\ref{4.12}) and (\ref{4.13}) finishes the proof. $\hfill\Box$

\medskip
We are now in position to give the proof of Theorem~1 as well as Corollary~1.

\medskip
{\bf Proof of Theorem~1.} It is sufficient to recall that eigenvalues of the problem ${\cal G}_k(q_2,q_3)$ coincide with zeros of the
characteristic function (\ref{2.12}), which has the form (\ref{4.1}) with
\begin{equation}\label{4.15}
g=(-1)^k\frac{\omega_2-\omega_3}4, \quad h= -\frac{\omega_2+\omega_3}4,
\end{equation}
and to apply Lemma~3. $\hfill\Box$

\medskip
{\bf Proof of Corollary~1.} Apply Corollary~2 and substitute (\ref{4.15}) into (\ref{4.6.2}) and (\ref{4.6.3}). Since the index ranges of the
subsequences in (\ref{4.6.1}) slightly differ from those of the subspectra in (\ref{1.17}), it remains to move one element from the first
subsequence to the second. $\hfill\Box$

\medskip
The next three lemmas will be used in Section~5 for proving Theorems~4 and~5.

\medskip
{\bf Lemma 4. }{\it The following representation holds:
\begin{equation}\label{4.16}
\Delta(\lambda)=-\frac\alpha{\lambda_0^0}\prod_{n=1}^\infty \frac{\lambda_n-\lambda}{\lambda_n^0},
\end{equation}
where
\begin{equation}\label{4.17}
\alpha=\lim_{\rho\to0}\frac{S(\rho)}{\rho^s}, \quad
\lambda_n^0=\left\{\begin{array}{cl} (z_n^0)^2, & z_n^0\ne0, \quad \\[3mm]
-1, & z_n^0=0, \end{array}\right. \quad n\ge0,
\end{equation}
and $s$ is the multiplicity of a zero of $S(\rho)$ at the origin, while $z_n^0$ are indexed as in Remark~1. }

\medskip
{\it Proof.} According to Theorem~5 in \cite{But22}, the function $\theta(\rho)$ of the form (\ref{4.2}) also has the representation
\begin{equation}\label{4.18}
\theta(\rho)=\alpha\exp(\beta\rho)\prod_{n\in A} \frac{z_n-\rho}{\mu_n}\exp\Big(\frac\rho{\mu_n}\Big),
\end{equation}
where $\alpha$ is defined in (\ref{4.17}), while $\mu_n$ are defined in (\ref{4.5}), and $\beta=s+\gamma,$ where
$$
\gamma=\lim_{\rho\to0}\frac{d}{d\rho}\ln\frac{S(\rho)}{\rho^s}=0
$$
since the function $S(\rho)$ and the number $s$ are even. Agreeing for definiteness that $z_{\pm n}^0=0$ for $n=\overline{0,s/2-1},$ and
taking the definition of $\mu_n$ into account, we rewrite (\ref{4.18}) in the form
\begin{equation}\label{4.19}
\theta(\rho)=\alpha\prod_{|n|=0}^{s/2-1}(\rho-z_n)\prod_{|n|=s/2}^\infty \frac{z_n-\rho}{z_n^0}\exp\Big(\frac\rho{z_n^0}\Big),
\end{equation}
where, in accordance with the definition of the index set $A,$ to each $|n|,$ including $n=0,$ there correspond two multipliers: $(\rho-z_n)$
and $(\rho-z_{-n}).$ Since, by virtue of the evenness of the functions $S(\rho)$ and $\theta(\rho),$ one can always assume that
$z_{-n}^0=-z_n^0$ and $z_{-n}=-z_n$ for $n\ge0,$ representation (\ref{4.19}) takes the form
\begin{equation}\label{4.20}
\theta(\rho)=\alpha\prod_{n=0}^{s/2-1}(\rho^2-z_n^2)\prod_{n=s/2}^\infty \frac{z_n^2-\rho^2}{(z_n^0)^2}.
\end{equation}
Recalling that $z_0=0$ and using the definitions of $\theta(\rho)$ and $\lambda_n^0,$ we arrive at (\ref{4.16}). $\hfill\Box$

\medskip
{\bf Lemma 5. }{\it Let $\{\lambda_n\}_{n\in{\mathbb N}}$ be an arbitrary complex sequence of the form (\ref{4.5.1}). Then the function
$\Delta(\lambda)$ constructed by formulae (\ref{4.16}) and (\ref{4.17}) has the form (\ref{4.1}).}

\medskip
{\it Proof.} Put $z_0:=0.$ Then the function $\theta(\rho):=\rho^2\Delta(\rho^2)$ has representation (\ref{4.20}). Then, after the
continuation $z_{-0}:=0$ and $z_{-n}:=-z_n$ for $n\in{\mathbb N},$ it takes the form (\ref{4.19}), which, in turn, is equivalent to
(\ref{4.18}) with $\beta=s.$ By virtue of the second part of Theorem~6 in \cite{But22}, the function $\theta(\rho)$ has the form
\begin{equation}\label{4.21}
\theta(\rho)=S(\rho)+P_0(\rho)S_0(\rho)+\int_{-3}^3 v(x)\exp(i\rho x)\,dx, \quad v(x)\in L_2(-3,3),
\end{equation}
where $P_0(\rho)$ is a polynomial of degree $0,$ i.e. $P_0(\rho)\equiv C - const,$ while $S_0(x)$ is a sine-type function determined by
(\ref{4.4.1}).

Further, since the functions $\theta(\rho)$ and $S(\rho)=P_1(\rho)S_0(\rho)$ are even, so is their difference
$f(\rho):=\theta(\rho)-S(\rho).$ Thus, we have
\begin{equation}\label{4.22}
C\frac{S_0(\rho)+S_0(-\rho)}2=f(\rho)-\int_{-3}^3 v(x)\cos\rho x\,dx, \;\; C\frac{S_0(-\rho)-S_0(\rho)}{2i}=\int_{-3}^3 v(x)\sin\rho x\,dx.
\end{equation}
Let $P_1(\rho)=h_1\rho+h_0.$ Then the evenness of $S(\rho)$ implies
\begin{equation}\label{4.23}
h_1\rho(S_0(\rho)+S_0(-\rho))=h_0(S_0(-\rho)-S_0(\rho)).
\end{equation}
Let $C\ne0.$ Then $f(\rho)$ is a sine-type function. On the other hand, substituting (\ref{4.22}) into (\ref{4.23}), we arrive at the
representation
$$
f(\rho)=\int_{-3}^3 v(x)\cos\rho x\,dx -\frac{ih_0}{h_1\rho}\int_{-3}^3 v(x)\sin\rho x\,dx.
$$
In particular, it gives $f(\rho)\to0$ as ${\rm Re}\rho\to\infty$ for any fixed ${\rm Im}\rho\in{\mathbb R},$ which contradicts the
left-hand estimate in (\ref{4.4.2}). This contradiction implies $C=0.$

Thus, the function $v(x)$ under the integral in (\ref{4.21}) is even. Hence, the function $\Delta(\lambda)$ has the form (\ref{4.1}) with
$w(x)=2v(x).$ $\hfill\Box$

\medskip
Finally, we give the uniform stability of recovering the function (\ref{4.1}) from its zeros.
For this purpose, along with $\Delta(\lambda),$ we consider
another function $\tilde\Delta(\lambda)$ of the same form (\ref{4.1}) and with the same coefficients $g$ and $h$ but with a different
function $\tilde w(x)$ under the integral. If a certain symbol~$\gamma$ denotes an object related to the function $\Delta(\lambda),$ then
this symbol with tilde $\tilde\gamma$ will denote the analogous object related to $\tilde\Delta(\lambda).$

The following lemma immediately follows from Theorem~7 in~\cite{But22}.

\medskip
{\bf Lemma 6. }{\it  For any $r>0,$ there exists $C_r>0$ such that the estimate
$$
\|\theta-\tilde\theta\|_{L_2(-\infty,\infty)}=\sqrt\pi\|w-\tilde w\|_{L_2(0,3)}\le C_r\|\{n(z_n-\tilde z_n)\}_{n\in{\mathbb N}}\|_{l_2}
$$
is fulfilled whenever $\|\{n(z_n-z_n^0)\}_{n\in{\mathbb N}}\|_{l_2}\le r$ and $\|\{n(\tilde z_n-z_n^0)\}_{n\in{\mathbb N}}\|_{l_2}\le r.$}
\\

{\large\bf 6. Solution of the inverse problems}
\\

Let us begin with Inverse Problem~2.

\medskip
Fix $k\in\{1,2\}$ and consider an arbitrary complex sequence $\{\mu_{n,k}\}_{|n|\in{\mathbb N}}$ of the form (\ref{1.19}). Supplement
it up to $\{\mu_{n,k}\}_{n\in{\mathbb Z}}$ with $\mu_{0,k}=0.$ Without loss of generality, we agree that
$$
\mu_{n,k}=\mu_{n+1,k}=\ldots=\mu_{n+m_{n,k}-1,k}, \quad n\in{\mathbb Z},
$$
where $m_{n,k}$ is the multiplicity of the value $\mu_{n,k}$ in the sequence $\{\mu_{n,k}\}_{n\in{\mathbb Z}}.$ Thus, the set
$$
{\cal S}_k:=\{n:\mu_{n,k}\ne \mu_{n-1,k},\, n\in{\mathbb Z}\}
$$
indexes elements of this sequence without account of multiplicity.

\medskip
{\bf Lemma 7. }{\it The functional sequence ${\cal C}_k:=\{c_{n,k}(x)\}_{n\in{\mathbb Z}}$ determined by the formula
\begin{equation}\label{5.3}
c_{l+\nu,k}(x):=\frac{d^\nu}{d\lambda^\nu}\cos\sqrt\lambda x\Big|_{\lambda=\mu_{l,k}}, \quad l\in{\cal S}_k, \quad \nu=\overline{0,m_{l,k}-1},
\end{equation}
is a Riesz basis in $L_2(0,2).$}

\medskip
{\it Proof.} According to Lemma~4 in \cite{Bond18}, the sequence ${\cal C}^\gamma:=\{\cos(\pi n+\gamma)x\}_{n\in{\mathbb Z}}$ is a Riesz basis
in $L_2(0,2)$ for any real $2\gamma/\pi\notin{\mathbb Z},$ which can be also extended to any non-real $\gamma$ as in Lemma~A1 in
\cite{ButHu21}. The asymptotics in (\ref{1.19}) implies $m_{n,k}=1$ for large $|n|$ and, hence,
$$
c_{n,k}(x)=\cos(\pi n+\sigma)x +O\Big(\frac1n\Big), \quad |n|\to\infty,
$$
uniformly in $x.$ Thus, the sequence ${\cal C}_k$ is quadratically close to the Riesz basis ${\cal C}^\sigma.$

It remains to show that ${\cal C}_k$ is complete in $L_2(0,2).$ For this purpose, consider the function
$$
B(\lambda)=\frac{G(\lambda)}{\Theta(\lambda)}, \quad G(\lambda):=\int_0^2g(x)\cos\sqrt\lambda x\,dx, \quad
\Theta(\lambda):=\lambda^{m_{0,k}}\prod_{\mu_{n,k}\ne0} \Big(1-\frac{\lambda}{\mu_{k,n}}\Big),
$$
where $g(x)\in L_2(0,2).$ As in Lemma~5 above, applying Theorem~6 in \cite{But22} one can show that the function $\Theta(\lambda)$
has the form
$$
\Theta(\lambda)=C(3\cos2\sqrt\lambda+1) +\int_0^2 V(x)\cos\sqrt\lambda x\,dx, \quad V(x)\in L_2(0,2), \quad C\ne0 - {\rm const},
$$
which means that $f(\rho):=\Theta(\rho^2)$ is a sine-type function of type $b=2.$ Hence, according to the left-hand inequality in
(\ref{4.4.2}), we obtain the estimate
$$
B(\rho^2)=o(1), \quad \rho\to\infty, \quad {\rm dist}(\rho,\{\pm\sqrt{\mu_{n,k}}\}_{n\in{\mathbb Z}})\ge\delta>0
$$
If $\{\mu_{n,k}\}_{n\in{\mathbb Z}}$ are zeros of the function $G(\lambda)$ with account of multiplicity, then the function
$B(\lambda)$ is entire. Thus, by virtue of the maximum modulus principle for analytic functions along with
Liouville's theorem, we arrive at $B(\lambda)\equiv0,$ which gives the required completeness. $\hfill\Box$

\medskip
For the proof of Theorems~2 and~3, we need the following auxiliary assertion.

\medskip
{\bf Lemma 8. }{\it Specification of two sequences $\{\mu_{n,1}\}_{|n|\in{\mathbb N}}$ and $\{\mu_{n,2}\}_{|n|\in{\mathbb N}}$ of the form
(\ref{1.19}) and (\ref{1.18.1}) uniquely determines the numbers $\omega_2$ and $\omega_3.$}

\medskip
{\it Proof.} According to (\ref{1.19}), we have $\mu_{n,k}=(\pi n +\sigma)^2 +2(-1)^n\gamma_k+\varkappa_n$ and, hence,
\begin{equation}\label{5.4}
\gamma_k=\frac12\lim_{n\to\infty}\Big((-1)^n(\mu_{n,k}-(\pi n +\sigma)^2)\Big), \quad k=1,2.
\end{equation}
For $k=1,2,$ transforming the second formula in (\ref{1.18.1}), we get
$$
\left[\begin{array}{c}\gamma_1\\ \gamma_2\end{array}\right]
 =\frac1{8\sqrt3}\left[\begin{array}{rr}-1&2\\2&-1\end{array}\right] \left[\begin{array}{c}\omega_2\\ \omega_3\end{array}\right].
$$
Solving this system, we calculate
\begin{equation}\label{5.5}
\left[\begin{array}{c}\omega_2\\ \omega_3\end{array}\right]
 =\frac8{\sqrt3}\left[\begin{array}{rr}1&2\\2&1\end{array}\right] \left[\begin{array}{c}\gamma_1\\ \gamma_2\end{array}\right],
\end{equation}
which finishes the proof. $\hfill\Box$

\medskip
{\bf Proof of Theorem~2.} Rewrite (\ref{2.12}) in the form
\begin{equation}\label{5.6}
\lambda\Delta_k(\lambda)=F_k(\lambda) +\int_0^2 w_k(x)\cos\sqrt\lambda x\,dx,
\end{equation}
where $F_k(\lambda)$ is defined in (\ref{1.16.2}).
For $k=1,2$ and $l\in{\cal S}_k,$ differentiate $\nu=\overline{0,m_{l,k}-1}$ times the relation (\ref{5.6}) and substitute
$\lambda=\mu_{l,k}$ into the obtained derivatives. Thus, we arrive at
\begin{equation}\label{5.8}
\beta_{n,k}=\int_0^2 w_k(x)c_{n,k}(x)\,dx, \quad n\in{\mathbb Z}, \quad k=1,2,
\end{equation}
where the functions $c_{n,k}(x)$ are defined in (\ref{5.3}), while
\begin{equation}\label{5.9}
\beta_{l+\nu,k}:=-F_k^{(\nu)}(\mu_{l,k}), \quad l\in{\cal S}_k, \quad \nu=\overline{0,m_{l,k}-1}.
\end{equation}
By virtue of Lemma~8, the function $F_k(\lambda)$ is uniquely determined by specifying both $\mu$-sub\-spectra. Therefore, so are the
sequences $\{\beta_{n,1}\}_{n\in{\mathbb Z}}$ and $\{\beta_{n,2}\}_{n\in{\mathbb Z}}.$

Thus, according to Lemma~7, relations (\ref{5.8}) uniquely determine both $w_k(x).$ Hence, by Lemma~2, the
potentials $q_2(x)$ and $q_3(x)$ are uniquely determined too. $\hfill\Box$

\medskip
{\bf Proof of Theorem~3.} Let the numbers $\omega_2$ and $\omega_3$ be determined by formulae (\ref{5.4}) and (\ref{5.5}) from the
given sequences $\{\mu_{n,1}\}_{n\in{\mathbb Z}}$ and $\{\mu_{n,2}\}_{n\in{\mathbb Z}}.$

Next, we show that the corresponding sequences $\{\beta_{n,1}\}_{n\in{\mathbb Z}}$ and $\{\beta_{n,2}\}_{n\in{\mathbb Z}}$
determined by (\ref{5.9}) belong to $l_2.$ Indeed, for large $|n|,$ we have $\beta_{n,k}=-F_k(\mu_{n,k})$ or, with account
of (\ref{1.16.2}),
\begin{equation}\label{5.10}
\beta_{n,k}=-\eta_{n,k}\frac{1+3\cos2\eta_{n,k}}2\sin\eta_{n,k} +\frac{\omega_2+\omega_3}4\cos2\eta_{n,k}-(-1)^k\frac{\omega_2-\omega_3}4,
\end{equation}
where, using (\ref{1.19}) along with the first relation in (\ref{1.18.1}), we calculate $\sin\sigma=\sqrt{2/3}$ and
\begin{equation}\label{5.12}
\cos2\eta_{n,k}=-\frac13-(-1)^n\frac{4\sqrt2}{3\pi n}\gamma_k+\frac{\varkappa_n}n, \quad
\frac{1+3\cos2\eta_{n,k}}2=-(-1)^n\frac{2\sqrt2}{\pi n}\gamma_k+\frac{\varkappa_n}n,
\end{equation}
\begin{equation}\label{5.14}
\sin\eta_{n,k}=(-1)^n\sqrt\frac23 +O\Big(\frac1n\Big), \quad |n|\to\infty.
\end{equation}
Substituting (\ref{1.19}), (\ref{5.12}) and (\ref{5.14}) into (\ref{5.10}), we arrive at
$$
\beta_{n,k}=\frac{4\gamma_k}{\sqrt3}-\frac{\omega_2+\omega_3}{12}-(-1)^k\frac{\omega_2-\omega_3}4+\varkappa_n,
$$
which along with the second relation in (\ref{1.18.1}) gives $\beta_{n,k}=\varkappa_n.$

Thus, by virtue of Lemma~7, there exists a unique pair of functions $w_1(x)$ and $w_2(x)$ obeying relations (\ref{5.8}).
Construct $q_2(x)$ and $q_3(x)$ using (\ref{2.14}) and consider the corresponding problems ${\cal G}_k(q_2,q_3),\;k=1,2.$ Since
$0\in\{\mu_{n,k}\}_{n\in{\mathbb Z}},$ formulae (\ref{1.16.2}), (\ref{5.3}), (\ref{5.8}) and (\ref{5.9}) yield
$$
\frac{\omega_2+\omega_3}4-(-1)^k\frac{\omega_2-\omega_3}4 =\int_0^2 w_k(x)\,dx,
$$
which along with (\ref{2.14}) implies (\ref{1.16.0}). Hence, due to Lemma~1, the function $\Delta_k(\lambda)$ determined by (\ref{2.12})
with the obtained $\omega_2,$ $\omega_3$ and $w_k(x)$ is the characteristic function of the problem ${\cal G}_k(q_2,q_3)$ for $k=1,2.$ On the
other hand, relations (\ref{5.3}), (\ref{5.8}) and (\ref{5.9}) mean that the sequences
$\{\mu_{n,1}\}_{|n|\in{\mathbb N}}$ and $\{\mu_{n,2}\}_{|n|\in{\mathbb N}}$ consist of zeros
(not all) of these $\Delta_1(\lambda)$ and $\Delta_2(\lambda),$ respectively. Thus, $\{\mu_{n,k}\}_{|n|\in{\mathbb N}}$ is a
$\mu$-subspectrum of the problem ${\cal G}_k(q_2,q_3)$ as $k=1,2.$ $\hfill\Box$

\medskip
This proof is constructive and gives the following algorithm for solving Inverse Problem~2.

\medskip
{\bf Algorithm 1.} Let the $\mu$-subspectra $\{\mu_{n,1}\}_{|n|\in{\mathbb N}}$ and $\{\mu_{n,2}\}_{|n|\in{\mathbb N}}$ be given. Then:
\begin{itemize}
\item[(i)] Using formulae (\ref{5.4}) and (\ref{5.5}), find the numbers $\omega_2$ and $\omega_3,$ which determine, in turn, the functions
$F_1(\lambda)$ and $F_2(\lambda)$ by formula (\ref{1.16.2});

\item[(ii)] For $k=1,2,$ put $\mu_{0,k}:=0$ and construct $\{c_{n,k}(x)\}_{n\in {\mathbb Z}}$ by (\ref{5.3}), and
$\{\beta_{n,k}\}_{n\in {\mathbb Z}}$ by (\ref{5.9});

\item[(iii)] Find the functions $w_1(x)$ and $w_2(x)$ by the formula
$$
w_k(x)=\sum_{n=-\infty}^\infty \beta_{n,k}c_{n,k}^*(x), \quad k=1,2,
$$
where $\{c_{n,k}^*(x)\}_{n\in {\mathbb Z}}$ is the biorthogonal basis to the basis $\{\overline{c_{n,k}(x)}\}_{n\in {\mathbb Z}};$

\item[(iv)] Construct the functions $q_2(x)$ and $q_3(x)$ by using (\ref{2.14}).
\end{itemize}

\bigskip
Now, we proceed with Inverse Problem~1.

\medskip
{\bf Proof of Theorem~4.} By necessity, the asymptotics (\ref{1.16.1}) was already established in Theorem~1. According to Lemma~4,
the characteristic functions $\Delta_1(\lambda)$ and $\Delta_2(\lambda)$ have the representation as in (\ref{1.19.1}). Thus, condition (ii)
easily follows from representation (\ref{2.12}).

For the sufficiency, construct the functions $\Delta_1(\lambda)$ and $\Delta_2(\lambda)$ by the second formula in~(\ref{1.19.1}) using the
given sequences $\{\lambda_{n,1}\}_{n\in{\mathbb Z}}$ and $\{\lambda_{n,2}\}_{n\in{\mathbb Z}}.$ By Lemma~5, these functions have the form
$$
\Delta_k(\lambda)=\Delta_0(\lambda)+(-1)^k\frac{\omega_2-\omega_3}{4\rho^2} -\frac{\omega_2+\omega_3}{4\rho^2}\cos2\rho +\int_0^3
w_k(x)\frac{\cos\rho x}{\rho^2}\,dx, \quad k=1,2,
$$
with some functions $w_1(x),w_2(x)\in L_2(0,3).$ Further, condition (ii) along with the Paley--Wiener theorem implies
$w_k(x)=0$ a.e.~on $(2,3)$ for $k=1,2,$ i.e. representation (\ref{2.12}) holds.

Construct $q_2(x)$ and $q_3(x)$ by formula (\ref{2.14}) and consider the corresponding problems ${\cal G}_k(q_2,q_3),$ $k=1,2.$ Then, as in
the proof of Theorem~3, one can show that (\ref{1.16.0}) holds. Hence, $\Delta_k(\lambda)$ is the characteristic function of the problem
${\cal G}_k(q_2,q_3)$ for $k=1,2.$ $\hfill\Box$

\medskip
{\bf Proof of Theorem~5.} Applying Lemma~6 to the characteristic functions (\ref{2.12}), we get
$$
\|w_k-\tilde w_k\|_{L_2(0,2)}\le C_r\|\{n(z_{n,k}-\tilde z_{n,k})\}_{n\in{\mathbb N}}\|_{l_2}, \quad k=1,2,
$$
as soon as the conditions of Theorem~5 are met. Hence, estimate (\ref{1.20}) follows from~(\ref{2.14}). $\hfill\Box$

\medskip
Finally, note that Algorithm~1 can be used also for solving Inverse Problem~1 because its input data are included into
those of Inverse Problem~2. However, it can be simplified in the following way because the complete spectra are known
and one can use an orthogonal basis.

\medskip
{\bf Algorithm 2.} Let the spectra $\{\lambda_{n,1}\}_{n\in{\mathbb N}}$ and $\{\lambda_{n,2}\}_{n\in{\mathbb N}}$ be given. Then:
\begin{itemize}
\item[(i)] Using formulae (\ref{1.16.2}), (\ref{5.4}) and (\ref{5.5}), calculate the functions
$F_1(\lambda)$ and $F_2(\lambda);$

\item[(ii)] Construct the functions $\Delta_1(\lambda)$ and $\Delta_2(\lambda)$by the second formula in~(\ref{1.19.1});

\item[(iii)] Find the functions $w_1(x)$ and $w_2(x)$ inverting the Fourier transform in (\ref{2.12}):
$$
w_k(x)=\sum_{n=0}^\infty \Big(\frac{\pi^2n^2}4\Delta_k\Big(\frac{\pi^2n^2}4\Big)-F_k\Big(\frac{\pi^2n^2}4\Big)\Big)\cos\frac{\pi nx}2, \quad
k=1,2;
$$

\item[(iv)] Construct the functions $q_2(x)$ and $q_3(x)$ by using (\ref{2.14}).
\end{itemize}

\bigskip

{\large\bf 7. An arbitrary compact tree with global delay}
\\

In this section, we extend the definition of the functional-differential operator (\ref{1.1}) to an arbitrary
compact tree (i.e. a graph without cycles). Although noncompact trees would not bring any essential additional difficulties,
we restrict ourself with the compact case.

Let ${\cal T}$ be a compact rooted tree,
with the set of vertices $V=\{v_0,v_1,\ldots,v_m\}$ and the set of edges $E=\{e_1,\ldots, e_m\}.$ The vertex
$v_0$ will be labelled as {\it root}. Denote by ${\rm deg}(v)$ the degree of $v\in V,$ i.e.
the number of edges incident on the vertex $v.$ Any vertex $v$ may be {\it boundary} one or {\it internal} one depending on
whether ${\rm deg}(v)=1$
or ${\rm deg}(v)>1,$ respectively. Without loss of generality, we assume that $v_0$ is a boundary vertex. Otherwise, one can split
the tree ${\cal T}$ into ${\rm deg}(v_0)$ subtrees and repeat our scheme for each of them.

For two geometrical points $t_1$ and $t_2$ on the tree ${\cal T},$ we write $t_1\le t_2$ if $t_1$ lies on the unique
simple path from the root $v_0$ to the point $t_2,$ whose length will be denoted by $|t_2|.$
In particular, we write $t_1<t_2$ if $t_1\le t_2$ and $t_1$ does not coincide with $t_2.$
In the latter case, we put $[t_1,t_2]:=\{z\in {\cal T}:\; t_1\le z\le t_2\}.$
If $e:=[v,w]$ is an edge, we call $v$ its {\it initial point}, $w$ its {\it end point} and say that $e$ {\it emanates} from $v$ and
{\it terminates} at $w.$ For any internal vertex~$v,$ we denote by $R(v)$ the set of edges emanating from $v,$
i.e. $R(v)=\{e\in E:\; e=[v,w],\, w\in V\}.$

Without loss of generality, we agree that vertices and edges are indexed so that $e_j=[v_{k_j},v_j]$ for $j=\overline{1,m},$
while $k_1=0$ and $\{v_1,\ldots,v_p\}$ are internal vertices, where $1\le p<m.$

The value $l=\max_{j=\overline{1,m}}|v_j|$ is called {\it height} of ${\cal T}.$ For $j=\overline{1,m},$
we denote by $l_j$ the length of the edge~$e_j,$ which is parameterized by $x\in[0,l_j]$ so that $x=0$ always corresponds to
the initial point of $e_j.$ For
$j=\overline{2,m},$ there exists a unique nonempty chain of edges
\begin{equation}\label{7.1}
e_{j_{1,j}}, \;\; e_{j_{2,j}}, \;\; \ldots \;\; e_{j_{\sigma_j,j}},
\end{equation}
that connects the edge $e_j$ with the root. In other words, the chain (\ref{7.1}) forms the path from $v_0$ to $v_{k_j}.$
Thus, we always have $j_{1,j}=1$ and $j_{\sigma_j,j}=k_j.$

While any function $y$ on ${\cal T}$ can be understood as the vector (\ref{3.1}), whose component $y_j(x)$
is defined on $e_j,$ we continue $y_j(x)$ for $j=\overline{2,m}$ to the chain (\ref{7.1}) by the formula
\begin{equation}\label{7.2}
y_j(x):=y_{j_{\nu,j}}(x+x_{\nu,j}), \;\; x\in [-x_{\nu,j}, -x_{\nu+1,j}), \;\; \nu=\overline{1,\sigma_j},
\quad x_{\nu,j}:=\sum_{k=\nu}^{\sigma_j}l_{j_{k,j}}, \;\; x_{\sigma_j+1,j}:=0.
\end{equation}

On edges of ${\cal T},$ we consider the functional-differential equations
\begin{equation}\label{7.3}
-y_j''(x)+q_j(x)y_j(x-a)=\lambda y_j(x), \quad 0<x<l_j, \quad j=\overline{1,m},
\end{equation}
where $a\in(0,l)$ is fixed and $q_j(x)\in L_2(0,l_j),$ while
\begin{equation}\label{7.3-1}
q_j(x)=0 \;\;{\rm a.e. \;\; on}\;\;(0,\min\{l_j,a-|v_{k_j}|\})\;\;{\rm if}\;\; |v_{k_j}|< a.
\end{equation}
Since $|v_{k_j}|=|v_j|-l_j,$ the latter is equivalent to $q_j=0$ if $|v_j|\le a,$ and $q_j(x)=0$ a.e. on $(0,a-|v_{k_j}|)$
if $|v_{k_j}|<a<|v_j|.$ All other $q_j(x)$ (i.e. when $|v_{k_j}|\ge a)$
are
arbitrary. We also assume that the unknown functions $y_j(x)$ in  (\ref{7.3}) for negative $x$ are determined by~(\ref{7.2}).

The assumptions  (\ref{7.2}) and (\ref{7.3-1}) make the system of equations (\ref{7.3}) well-defined. Indeed, according to (\ref{7.2}),
the function $y_j(x)$ is additionally defined for $x\in[-|v_{k_j}|,0)$ as soon as $j=\overline{2,m}.$ On the other hand, (\ref{7.3-1})
actually eliminates $y_j(x-a)$ in (\ref{7.3}) if $x-a<-|v_{k_j}|.$

Along with (\ref{7.3}), we consider the standard matching conditions at the internal vertices
\begin{equation}\label{7.4}
y_j(1)=y_k(0) \;\; {\rm for \;\, all}\;\; e_k\in R(v_j),
\quad y'_j(1)=\sum_{e_k\in R(v_j)} y_k'(0), \quad j=\overline{1,p},
\end{equation}
as well as the Dirichlet or the Neumann boundary conditions at the boundary ones
\begin{equation}\label{7.5}
y_1^{(\nu_1)}(0)=0, \quad y_j^{(\nu_j)}(l_j)=0, \quad j=\overline{p+1,m},
\end{equation}
where $\nu_j\in\{0,1\}.$ Denote by ${\cal G}$ the boundary value problem on ${\cal T}$ consisting of (\ref{7.2})--(\ref{7.5}).

Finally, let us give some illustrative examples.

\medskip
{\bf Example 1.} Let $l_j=1,\;j=\overline{1,m},$ $a\in(0,2].$ Then relations (\ref{7.2}) can be replaced with
$$
y_j(x-a)=\left\{\begin{array}{cc}
y_{k_j}(x-a+1), & \max\{0,a-1\}<x<\min\{a,1\}, \quad  \\[3mm]
y_{k_{k_j}}(x-a+2), & 0<x<\max\{0,a-1\}, \quad |v_j|>2,
\end{array}\right. \quad j=\overline{2,m},
$$
while (\ref{7.3-1}) implies $q_1(x)=0$ a.e. on $(0,\min\{a,1\})$ and $q_j(x)=0$ a.e. on $(0,\max\{0,a-1\})$ whenever $|v_j|=2$
(i.e. $k_j=1).$ In particular, the problem ${\cal G}$ coincides with the problem (\ref{1.8})--(\ref{1.11}) if $p=1.$

\medskip
In what follows, we consider the simplest non-star tree, i.e. when $m=5$ and $p=2$ (see Fig.~4).
For clarity, we will also assume that $l_j=1$ for all $j.$
\begin{center}
\unitlength=0.7mm
\begin{picture}(180,97)
 \put(127,64.8){\line(1,1){27.7}}
 \put(127,64.8){\line(2,-1){50}}
 \put(80,46){\line(-1,0){50}}
 \put(80,46){\line(5,2){46.65}}
 \put(80,46){\line(2,-3){23}}

 \put(30,46){\circle*{1}}
 \put(80,46){\circle*{1}}
 \put(155.1,92.9){\circle*{1}}
 \put(127,64.8){\circle*{1}}
 \put(103,11.5){\circle*{1}}
 \put(177,39.8){\circle*{1}}

 \put (51,41){\small $e_1$}
 \put (104,52.5){\small $e_2$}
 \put (85,25){\small $e_3$}
 \put (142,77){\small $e_4$}
 \put (148,47){\small $e_5$}

 \put (22,44){\small $v_0$}
 \put (76,50){\small $v_1$}
 \put (121,68){\small $v_2$}
 \put (102.5,6.5){\small $v_3$}
 \put (155,94){\small $v_4$}
 \put (178,37){\small $v_5$}

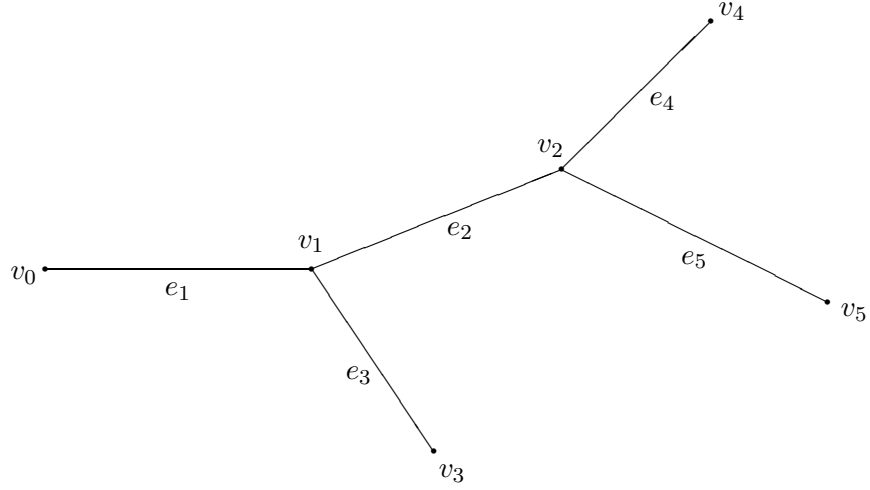
 \put (0,0){\small Fig. 4. The simplest non-star ${\cal T}$}

\end{picture}
\end{center}

{\bf Example 2.} Let $a=1.$ Then the problem ${\cal G}$ takes the form
$$
\;\;\; -y_1''(x)=\lambda y_1(x), \quad 0<x<1,
$$
$$
-y_j''(x)+q_j(x)y_1(x)=\lambda y_j(x), \quad 0<x<1, \quad j=2,3,
$$
$$
-y_j''(x)+q_j(x)y_2(x)=\lambda y_j(x), \quad 0<x<1, \quad j=4,5.
$$
$$
y_1(1)=y_2(0)=y_3(0), \quad y_1'(1)=y_2'(0)+y_3'(0),
$$
$$
y_2(1)=y_4(0)=y_5(0), \quad y_2'(1)=y_4'(0)+y_5'(0),
$$
$$
y_j^{(\nu_j)}(1-\delta_{1,j})=0, \quad  j=1,3,4,5,
$$
where, in particular, the first three lines are induced by relations (\ref{7.2}) and (\ref{7.3}).

\medskip
{\bf Example 3.} For $a=2,$ relations (\ref{7.2}) and (\ref{7.3}) can be replaced with
$$
-y_j''(x)=\lambda y_j(x), \;\; j=\overline{1,3}, \quad
-y_j''(x)+q_j(x)y_1(x)=\lambda y_j(x), \;\; j=4,5, \qquad 0<x<1.
$$

\medskip
{\bf Example 4.} For $a=3/2,$ relations (\ref{7.2}) and (\ref{7.3}) imply
$$
-y_1''(x)=\lambda y_1(x), \quad 0<x<1,
$$
$$
\left.\begin{array}{rc}
\displaystyle -y_j''(x)=\lambda y_j(x), & \displaystyle 0<x<\frac12, \\[3mm]
\displaystyle -y_j''(x)+q_j(x)y_1\Big(x-\frac12\Big)=\lambda y_j(x), & \displaystyle \frac12<x<1,
\end{array}\right\} \quad j=2,3, \qquad
$$
$$
\left.\begin{array}{cc}
\displaystyle -y_j''(x)+q_j(x)y_1\Big(x+\frac12\Big)=\lambda y_j(x), & \displaystyle 0<x<\frac12, \\[3mm]
\displaystyle -y_j''(x)+q_j(x)y_2\Big(x-\frac12\Big)=\lambda y_j(x), & \displaystyle \frac12<x<1,
\end{array}\right\} \quad j=4,5, \qquad
$$
under the natural continuity conditions
$$
y_j^{(\nu)}\Big(\frac12-0\Big) =y_j^{(\nu)}\Big(\frac12+0\Big), \quad \nu=0,1, \quad j=\overline{2,5}.
$$
\\

{\bf Funding.} This research was supported by Russian Science Foundation, Grant No. 22-21-00509, https://rscf.ru/project/22-21-00509/
\\

{\bf Acknowledgement.} The author is grateful to Maria Kuznetsova, who has carefully read the manuscript and made valuable comments.

\end{document}